\documentclass{article}\begin{document}\centerline{\bf\ Jacobian Elliptic Functions, Continued Fractions and Ramanujan Quantities.}\vskip .10in

\centerline{\bf Nikos Bagis}

\centerline{Department of Informatics}

\centerline{Aristotle University of Thessaloniki}
\centerline{54124 Thessaloniki, Greece}\centerline{ and}

\centerline{\bf M.L. Glasser}
\centerline{Department of Physics}
\centerline{ Clarkson University}
\centerline{Potsdam, New York (USA)}

\begin{quote}
\begin{abstract}
In this article we present ways to evaluate certain sums, products and continued fractions using tools from the theory of elliptic functions.    The specific results appear to be new, although similar ones can be found in the literature; in most cases the methods applied are different.
\end{abstract}

\bf keywords \rm{Jacobian Elliptic Functions; Continued Fractions; Ramanujan Quantities}

\end{quote}

\newpage
\section{Introduction}
\label{intro}
Let
\begin{equation}
\left(a;q\right)_k=\prod^{k-1}_{n=0}(1-aq^n)
\end{equation}
Then we define
\begin{equation}
f(-q)=(q;q)_\infty 
\end{equation}
and
\begin{equation}
\phi(-q)=(-q;q)_\infty
\end{equation}
Also let
\begin{equation}
K(x)=\int^{\pi/2}_{0} \frac{1}{\sqrt{1-x^2\sin^2(t)}}dt
\end{equation}
be the elliptic integral of the first kind see also [9],[7]\\
The ratio
\[
\prod^{\infty}_{n=1}(1-q^n)=\frac{\prod^{\infty}_{n=1}(1-q^{2n})}{\prod^{\infty}_{n=1}(1+q^n)}
\]
is known. For example from [9]:
\begin{equation}
\prod^{\infty}_{n=1}(1-q^{2n})^6=\frac{2kk^{'}K(k)^3}{\pi^3q^{1/2}}
\end{equation}
and from
\begin{equation}
q^{1/3}\prod^{\infty}_{n=1}(1+q^n)^8=2^{-4/3}\left(\frac{k}{1-k^2}\right)^{2/3} 
\end{equation}
we get
\begin{equation}
\prod^{\infty}_{n=1}(1-q^n)^8=\frac{2^{8/3}}{\pi^4}q^{-1/3}k^{2/3}(k^{'})^{8/3}K(k)^4 
\end{equation}
The variable $k$ is defined from the equation  
\begin{equation}
\frac{K(k^{'})}{K(k)}=\sqrt{r}
\end{equation}
where $r$ is positive , $q=e^{-\pi \sqrt{r}}$ and $k^{'}=\sqrt{1-k^2}$. Note also that whenever $r$ is positive rational, the $k$ are algebraic numbers.\\
Some examples  of such products are
\[
\prod^{\infty}_{n=1}(1+e^{-n\pi \sqrt{2/5}})^8=1/8(7+3\sqrt{5}) e^{\pi/3 \sqrt{2/5}}
\]
in the same way for $r=3$
\[
\prod^{\infty}_{n=1}(1+e^{-n\pi \sqrt{3}})^8=\frac{e^{\pi/\sqrt{3}}}{2^{2/3}(26+15\sqrt{3})^{1/3}}
\]  
To evaluate $\prod^{\infty}_{n=1}(1-e^{-n\pi \sqrt{3}})^8$ we find from tables that 
\[
k^{'}=\frac{\sqrt{2+\sqrt{3}}}{2}, K(k)=\frac{3^{1/4}\Gamma^3(1/3)}{2^{7/3}\pi}
\]  
and the result is
\[
\prod^{\infty}_{n=1}(1-e^{-n\pi \sqrt{3}})^8=\frac{3(2+\sqrt{3})e^{\pi/\sqrt{3}} \Gamma^{12}(1/3)}{1024\pi^8}
\]
where $\Gamma$ is the Euler`s Gamma Function.

\section{Theta Functions}
\label{intro}
\textbf{Lemma 1.}
Let $k$ be the root of $K(k^{'})/K(k)=\sqrt{r}$, $q=e^{-\pi\sqrt{r}}$ then:
\begin{equation}
k=\frac{8q^{1/2}\phi(-q)^{12}}{1+\sqrt{1+64q\phi(-q)^{24}}}
\end{equation}
\begin{equation}
K(k)=\frac{f(-q)^2\pi\sqrt{1+\sqrt{1+64q\phi(-q)^{24}}}}{2\sqrt{2}\phi(-q)^2}
\end{equation}
\\
\textbf{Proof.}
From the relations (2),(3),(4) and equations (5),(6),(7), we have 
\[k^{'}=\frac{\pi f(-q)^2}{2 K(k)\phi(-q)^2}\]
and
\[
k=\frac{q^{1/2} \pi^2 f(-q)^4 \phi(-q)^8}{K(k)^2}
\]
using the relation $k^{'}=\sqrt{1-k^2}$ we get the result.
\\
\textbf{Theorem 1.} If $r$ is  positive as in the introduction, set $k_{11}=k$, $k_{12}=\sqrt{1-k_{11}^2}$, $k_{21}=(2-k_{11}^2-2k_{12})/k_{11}^2$, $k_{22}=\sqrt{1-k_{21}^2}$, then 
\begin{equation}
\sum^{\infty}_{n=-\infty}q^{n^2+2mn}=2^{1/6}q^{-m^2}\frac{(k_{11}k_{22})^{1/3}}{(k_{12}k_{21})^{1/6}}\sqrt{\frac{K(k_{11})}{\pi}}
\end{equation}
\begin{equation}
\sum^{\infty}_{n=-\infty}q^{n^2+(2m+1)n}=2^{5/6}q^{-1/4(2m+1)^2}\frac{(k_{11}{k_{12}}{k_{21}})^{1/6}}{k_{22}^{1/3}}\sqrt{\frac{K(k_{11})}{\pi}} 
\end{equation} 
where 
\[
k=\frac{8q^{1/2}\phi(-q)^{12}}{1+\sqrt{1+64q\phi(-q)^{24}}}
\]
and $m\in\bf Z$.
\\
\textbf{Proof}. From [2]:
\begin{equation}
S_z:=\sum^{\infty}_{n=-\infty}q^{n^2+zn}=\prod^{\infty}_{n=0}(1-q^{2n+2})(1+q^{2n+1+z})(1+q^{2n+1-z})
\end{equation}
and (5), one has, after some rearrangement,
\begin{equation}
S_{2m}=(2kk^{'})^{1/6}q^{-1/12}\sqrt{K/\pi}\prod^{m-1}_{n=0}\left(\frac{1+q^{2(n-m)+1}}{1+q^{2n+1}}\right)\prod^{\infty}_{n=0}(1+q^{2n+1})^2
\end{equation} 
But  the last product in (14) is
\begin{equation}
\prod^{\infty}_{n=1}\frac{(1+q^n)^2}{(1+q^{2n})^2}=q^{1/12}\frac{k_{11}^{1/6}k_{22}^{1/3}}{k_{21}^{1/6}k_{12}^{1/3}}
\end{equation}
by (6) and the Landen transformation. This proves (11). \\The proof of (12) proceeds similarly. 
\\
\textbf{Application.}
When $a \in \bf Z$ and $q$ is an arbitrary number such that $|q|<1$, then
\\
$\frac{1}{1-}\frac{q^{2a+2}}{1+}\frac{q^{2a+2}(1-q^2)}{1-}\frac{q^{2a+6}}{1+}\frac{q^{2a+4}(1-q^4)}{1-}...+$
\begin{center}
$+\frac{q^{-2a}}{1-}\frac{q^{-2a+2}}{1+}\frac{q^{-2a+2}(1-q^2)}{1-}\frac{q^{-2a+6}}{1+}\frac{q^{-2a+4}(1-q^4)}{1-}...=$
\end{center}
\begin{flushright}
$=2^{5/6}q^{-1/4(2a+1)^2}\frac{(k_{11}k_{12}k_{21})^{1/6}}{k_{22}^{1/3}}\sqrt{K(k_{11})/\pi}$
\end{flushright}
\textbf{Proof.} From [8] pg 596 we have 
\begin{equation}
\sum^{\infty}_{k=0}(-c)^k q^{k(k+1)/2}=\frac{1}{1+}\frac{cq}{1+}\frac{c(q^2-q)}{1+}\frac{cq^3}{1+}\frac{c(q^4-q^2)}{1+}...
\end{equation}
We set
\begin{equation}
M(c,q):=\frac{1}{1-}\frac{cq}{1+}\frac{c(q-q^2)}{1-}\frac{cq^3}{1+}\frac{c(q^2-q^4)}{1-}...
\end{equation}
then 
\[
M(c,q)+1/cM(1/c,q)=\sum^{\infty}_{k=-\infty}c^kq^{k(k+1)/2}
\]
If now $q=e^{-\pi n \sqrt{r}}$, $r-$positive real, set $c=q^a,a\in \textbf{Z}$ then 
\[
M(q^a,q)+q^{-a} M(q^{-a},q)=\sum^{\infty}_{k=-\infty}q^{k^2/2+(a+1/2)k}
\]
\[
M(q^{2a},q^2)+q^{-2a} M(q^{-2a},q^2)=\sum^{\infty}_{k=-\infty}q^{k^2+(2a+1)k}
\]
The application follows from Theorem 1.
\\
Note. One can evaluate $M(c,q)$ exactly when $c=-q^a$, $a$-odd integer (the proof is easy):
\begin{flushleft}
$\frac{q^{(a+1)^2/4}}{1-}\frac{q^{a+2}}{1-}\frac{q^a(q^4-q^2)}{1-}\frac{q^{a+6}}{1-}\frac{q^a(q^8-q^4)}{1-}\frac{q^{a+10}}{1-}\frac{q^a(q^{12}-q^6)}{1-}\ldots$=
\end{flushleft}
\begin{flushright}
$=1/2-\sum^{(a-1)/2}_{k=0}q^{k^2}+\frac{\vartheta_3(q)}{2}$
\end{flushright}
where $\vartheta_3(q)=\sum^{\infty}_{k=-\infty}q^{k^2}=\sqrt{\frac{2K(k_{11})}{\pi}}$
\\
\textbf{Lemma 2.}
\begin{equation}
\sum^{\infty}_{k=1}\frac{\cosh(2tk)}{k\sinh(\pi ak)}=\log(P_0)-\log(\vartheta_4(it,e^{-a\pi}))
\end{equation}
Where $P_0=\prod^{\infty}_{n=1}(1-e^{-2n\pi a})$ and $\vartheta_4(u,q)=1+\sum^{\infty}_{n=1}(-1)^nq^{n^2}\cos(2nu)$
\\
\textbf{Proof.} From [2] pg.170 relation (13-2-12) and the definition of theta functions we have
\begin{equation}
\vartheta_4(z,q)=\prod^{\infty}_{n=0}(1-q^{2n+2})(1-q^{2n-1}e^{2iz})(1-q^{2n-1}e^{-2iz})
\end{equation}
By taking the logarithm of both sides and expanding the logarithm of the individual terms in a power series it is simple to show (18) from (19), where $q=e^{-\pi a}$, $a$ positive real.\\
Now it is well known [7] that if $|q|<1$
\begin{equation}
R(q)=\frac{1}{1+}\frac{q}{1+}\frac{q^2}{1+}\frac{q^3}{1+}\cdots  
\end{equation}
satisfies the famous Roger`s Ramanujan identity:
\begin{equation}
R^{*}(q)=q^{-1/5}R(q)=\frac{(q;q^5)_{\infty}(q^4;q^5)_{\infty}}{(q^2;q^5)_{\infty}(q^3;q^5)_{\infty}}=\prod^{\infty}_{n=1}(1-q^n)^{X_2(n)}
\end{equation}
where $X_2(n)$ is the Legendre symbol $\left(\frac{n}{5}\right)$.\\
\textbf{Proposition 1.}
\begin{equation}
R(e^{-x})=e^{-x/5}\frac{\vartheta_4(3ix/4,e^{-5x/2})}{\vartheta_4(ix/4,e^{-5x/2})},  x>0
\end{equation}
\textbf{Proof.} From (21) taking logarithms in both sides we get 
\begin{equation}
R(e^{-x})=\exp\left(-x/5-\sum^{\infty}_{n=1}\frac{1}{n}\frac{e^{4nx}-e^{3nx}-e^{2nx}+e^{nx}}{e^{5nx}-1}\right), x>0
\end{equation}
which can be written as
\begin{equation}
R(e^{-x})=e^{-x/5}\frac{\exp\left(\sum^{\infty}_{n=1}\frac{\cosh(nx/2)}{n \sinh(5nx/2)}\right)}{\exp\left(\sum^{\infty}_{n=1}\frac{\cosh(3nx/2)}{n \sinh(5nx/2)}\right)}
\end{equation}
finally using Lemma 2 we get the result.
\[
\]
In the same way as above we have that if
\begin{equation}
H(q)=\frac{q^{1/2}}{(1+q)+}\frac{q^2}{(1+q^3)+}\frac{q^4}{(1+q^5)+}\frac{q^6}{(1+q^7)+}\cdots
\end{equation} 
then
\begin{equation}
H(e^{-x})=\exp\left(-x/2-\sum^{\infty}_{n=1}\frac{1}{n}\frac{e^{7nx}-e^{5nx}-e^{3nx}+e^{nx}}{e^{8nx}-1}\right), x>0
\end{equation}
\begin{equation}
H(e^{-x})=e^{-x/2}\frac{\vartheta_4(3ix/2,e^{-4x})}{\vartheta_4(ix/4,e^{-4x})},  x>0
\end{equation}
In the same way as above we have
\[
\] 
\textbf{Definition 1.} In general if $q=e^{-\pi \sqrt{r}}$ where $a,p,r>0$ we denote `Agile` the quantity 
\begin{equation}
[a,p;q]=(q^{p-a};q^p)_{\infty}(q^a;q^p)_{\infty}
\end{equation}
\[
\]
\textbf{Observation 1.(Unproved)}\\
If $q=e^{-\pi \sqrt{r}}$, $a,b,p,r$ positive rationals then
\begin{equation}
q^{p/12-a/2+a^2/{(2p)}}[a,p;q]=Algebraic
\end{equation}
\[
\]
\textbf{Definition 2.}\\We call  
\begin{equation}
R^{*}(a,b,p;q)=\frac{[a,p;q]}{[b,p;q]}
\end{equation}
\\`Ramanujan's Quantity`
 because many of Ramanujan's continued fractions can be put in this form
\[
\]
\textbf{Theorem 2.}\\ 
If $a,b,p,r$ are positive rationals, then
\begin{equation}
R(a,b,p;q):=q^{-(a-b)/2+(a^2-b^2)/(2p)}R^{*}(a,b,p;q)=Algebraic
\end{equation} 
\[
\]
\textbf{Proof.}
Eq.(31) follows easy from the Observation 1 and the definitions 1,2. 
\[
\]
One example is the Rogers-Ramanujan continued fraction
\begin{equation}
q^{1/5}R^*(1,2,5;q)=R^{*}(q)q^{1/5}=R(q)
\end{equation}
\[
\]
\textbf{Theorem 3.} For all positive reals $a,b,p,x$ 
\begin{equation} R(a,b,p;e^{-x})=\exp\left(-x\frac{a^2-b^2}{2p}+x\frac{a-b}{2}\right)\frac{\vartheta_4((p-2a)ix/4,e^{-px/2})}{\vartheta_4((p-2b)ix/4,e^{-px/2})}=
\end{equation}
\begin{equation}
=\exp\left(-x(\frac{a^2-b^2}{2p}-\frac{a-b}{2})-\sum^{\infty}_{n=1}\frac{1}{n}\frac{e^{anx}+e^{(p-a)nx}-e^{(p-b)nx}-e^{bnx}}{e^{pnx}-1}\right)
\end{equation}
\\
\textbf{Proof.} We use Lemma 2 to rewrite $R$ in the form 
\[
R(a,b,p;e^{-x})=\exp\left(-x\frac{a^2-b^2}{2p}+x\frac{a-b}{2}\right)\frac{\exp\left(\sum^{\infty}_{n=1}\frac{\cosh(nx(p-2b)/2)}{n \sinh(pnx/2)}\right)}{\exp\left(\sum^{\infty}_{n=1}\frac{\cosh((p-2a)nx/2)}{n \sinh(pnx/2)}\right)}
\]
\textbf{Theorem 4.} If $a,b,p$ are positive integers, then
\begin{equation}
R^{*}(a,b,p;q)=\prod^{\infty}_{n=1}(1-q^n)^{X_2(n)}
\end{equation}
where
\begin{equation}
X_2(n)= \left\{\begin{array}{cc}
                                    1, \mbox{  }  n\equiv(p-a)modp \\
                                   -1, \mbox{  }  n\equiv(p-b)modp \\
                                    1, \mbox{  }  n\equiv a mod p \\
                                   -1, \mbox{  }  n\equiv bmodp \\
                                    0, \mbox{  }  p|n 
\end{array}\right\}
\end{equation}
\textbf{Proof.}\\ 
Use Theorem 3. Take the logarithms and expand the product (35) the proof is easy.\\
\textbf{Applications.}\\
\textbf{1)} In general one has the identity\\
\begin{equation}      
[a,p;q]=\frac{1}{f(-q^p)}\sum^{\infty}_{n=-\infty}(-1)^n q^{pn^2/2+(p-2a)n/2}
\end{equation}
\\
Hence if $a,b,p$ are positive integers then 
\[
\frac{q^{p/12-a/2+a^2/(2p)}}{f(-q^p)}\sum^{\infty}_{n=-\infty}(-1)^nq^{pn^2/2+(p-2a)n/2}=Algebraic
\]
\\ 
\textbf{2)} Let $M(c,q)=\sum^{\infty}_{n=0}c^nq^{(n+1)n/2}$ then
\begin{equation}
M(-q^{-a},q^p)-q^aM(-q^a,q^p)=f(-q^p)[a,p;q]
\end{equation}
Hence
\begin{equation} 
\frac{M(-q^{-a},q^p)-q^a M(-q^a,q^p)}{M(-q^{-b},q^p)-q^b M(-q^b,q^p)}=R^{*}(a,b,p;q)
\end{equation} 
Thus every Ramanujan quantity can calculated by the function\\$\sum^{\infty}_{n=0}c^nq^{(n+1)n/2}$.
For example if $a=1,b=2,p=5$ we have
\begin{equation}
\frac{M(-q^{-1},q^5)-q M(-q,q^5)}{M(-q^{-2},q^5)-q^2M(-q^2,q^5)}=\frac{1}{1+}\frac{q}{1+}\frac{q^2}{1+}\cdots
\end{equation}
\[
\]
\textbf{3)} Let
\begin{equation} \tau_0(a,q)=\sqrt{\frac{\pi}{K(k_1)}}f(-q)q^{a^2/2-a/2+1/8}
\frac{[2a,2;q]}{[a,1;q]}
\end{equation}
then observe numerically that if $|q|<1$ and $n$ is integer:
\begin{equation}
\tau_0(a,q)=\tau_0(n\pm a;q)
\end{equation}
\begin{equation}
\left(\frac{\partial \tau_0(a,q)}{\partial a}\right)_{a \in Z}=0 
\end{equation}
for every $|q|<1$.
\[
\]
More generally if we define
\begin{equation}
\tau^*(a,p;q)=\sqrt{\frac{\pi}{K(k_1)}}f(-q^p)q^{a^2/(2p)-a/2+p/8} \frac{[2a,2p;q]}{[a,p;q]}
\end{equation}
Then if $a,p\in \bf R$, $n\in \bf N$
\begin{equation}
\tau^*(a,p;q)=\tau^*(np\pm a,p;q)
\end{equation}
Where
\begin{equation}
\psi^*(a,p;q)=\sum^{\infty}_{n=-\infty}q^{pn^2/2+(p-2a)n/2}
\end{equation}
\begin{equation}
\psi^*(a,p;q)=f(-q^p)(-q^a;q^p)_{\infty}(-q^{p-a};q^p)_\infty=f(-q^p) \frac{[2a,2p;q]}{[a,p;q]}
\end{equation}
For the $\tau$ function we have  next:
\[
\]
\textbf{Theorem 5.} When $a,b$ are real and $n$ integer
\begin{equation}
\tau^*(a,\frac{|a\pm b|}{n};q)=\tau^*(b,\frac{|a\pm b|}{n};q)
\end{equation}
When $a,b\in \bf Q$ we have
\begin{equation}
\tau^*(1/a,\frac{gcd(a,b)}{ab};q)=\tau^*(1/b,\frac{gcd(a,b)}{ab};q)
\end{equation}
Also if $a,b,r \in \bf Q$ then $\tau^*(a,p;e^{-\pi \sqrt{r}})=$Algebraic  
\[
\]
\textbf{Proof.}
We only prove (48). Rewrite (44) in the form
\begin{equation} 
q^{a^2/(2p)-a/2}\psi^*(a,p;q)=q^{b^2/(2p)-b/2}\psi^*(b,p;q)
\end{equation} 
where $p=|a\pm b|/(2m),m$ positive integer.\\ \
We rewrite (48) in the form
\begin{equation} q^{(B_1^2-A^2)/(4A)}\sum^{\infty}_{n=-\infty}q^{n^2A+B_1n}=q^{(B_2^2-A^2)/(4A)}\sum^{\infty}_{n=-\infty}q^{n^2A+B_2n}
\end{equation}
where $A=p/2, B_2=\pm B_1+ 2m A$ and then finally we arrive to an identity for $\vartheta_3(z,q)$ functions. The result follows easily from the periodicity of $\vartheta_3(z,q)$ with respect to $z$ and the fact that $m$ is an integer.
\[
\]
\textbf{The first order derivatives of some of Ramanujan`s continued fractions.}
\\Observe that if
\begin{equation}
R_1(q)=\frac{q^{1/5}}{1+}\frac{q}{1+}\frac{q^2}{1+}\frac{q^3}{1+}\cdots=q^{1/5}\frac{(q;q^5)_{\infty}(q^4;q^5)_{\infty}}{(q^2;q^5)_{\infty}(q^3;q^5)_{\infty}}
\end{equation}    
\begin{equation}
R_2(q)=\frac{q^{1/3}}{1+}\frac{q+q^2}{1+}\frac{q^2+q^4}{1+}\frac{q^3+q^6}{1+}\cdots=q^{1/3}\frac{(q;q^6)_{\infty}(q^5;q^6)_{\infty}}{(q^3;q^6)_\infty^2}
\end{equation}    
\begin{equation}
R_3(q)=\frac{q^{1/2}}{(1+q)+}\frac{q^2}{(1+q^3)+}\frac{q^4}{(1+q^7)+}\cdots=q^{1/2}\frac{(q;q^8)_{\infty}(q^7;q^8)_{\infty}}{(q^3;q^8)_\infty (q^5;q^8)_{\infty}}
\end{equation}    
then all these have derivatives
\begin{equation} 
R_{1,2,3}^{'}(q)\frac{q\pi^2}{K(k_r)^2}=Algebraic
\end{equation}
whenever $q=e^{-\pi \sqrt{r}}$ and $r$ is a positive rational.
\[
\]
\textbf{Observation 2.} If $a,b,p,r$ are positive rationals then
\begin{equation}
\frac{d}{dq}R(a,b,p;q)=\frac{K(k_r)^2}{q \pi^2}Algebraic
\end{equation}
\begin{equation}
\frac{d}{dq}\left(q^{p/12-a/2+a^2/(2p)}[a,p;q]\right)=\frac{K(k_r)^2}{q \pi^2}Algebraic
\end{equation}
\[
\]
\textbf{Examples.}
\begin{equation}
\left(\frac{d}{dq}R(1,2,4;q)\right)_{q=e^{-\pi}}=\frac{e^{\pi}\Gamma(1/4)^4}{64\cdot2^{5/8}\pi^3}
\end{equation}
\[
\left(\frac{d}{dq}R(1,2,5;q)\right)_{q=e^{-\pi}}=\frac{e^{\pi}\Gamma(1/4)^4}{16\pi^3}\rho
\]
where $\rho$ is root of 
\[
16-240t^2+800t^3-2900t^4-6000t^5-6500t^6+17500t^7+625t^8=0
\]
From [7] pg.24 we have the next
\[
\]
\textbf{Proposition.} Suppose that $a,b$ and $q$ are complex numbers with $\left|a b \right|<1$ and $\left|q\right|<1$ or that $a=b^{2m+1}$ for some integer $m$. Then  
\[
P(a,b,q):=\frac{(a^2 q^3;q^4)_{\infty} (b^2q^3;q^4)_{\infty}}{(a^2q;q^4)_{\infty}(b^2q;q^4)_{\infty}}=
\]
\[
\frac{1}{(1-ab)+}\frac{(a-bq)(b-aq)}{(1-ab)(q^2+1)+}\frac{(a-bq^3)(b-aq^3)}{(1-ab)(q^4+1)+...}
\]
One can easily see that
\begin{equation}
P(q^A,q^B,q^{A+B})=\frac{(q^a;q^p)_{\infty}(q^{2p-a};q^p)_{\infty}}{[b,p;q]} 
\end{equation}
where $a=2A+3p/4$, $b=2B+p/4$ and $p=4(A+B)$.
\\
Let now define 
\begin{equation}
_2\phi_1[a,b;c;q,z]:=\sum^{\infty}_{n=0} \frac{(a;q)_n (b;q)_n}{(c;q)_n}\frac{z^n}{(q;q)_n}
\end{equation}
and
\begin{equation}
\psi(a,q,z):=\sum^{\infty}_{n=0}\frac{(a;q)_n}{(q,q)_n}z^n=_2\phi_1[a,0,0,q,z]
\end{equation}
Then  
\[
\] 
\textbf{Theorem 6.}\\
i) If $a=2A+3p/4$, $b=2B+p/4$ and $p=4(A+B)$, $|q|<1$
\begin{equation}
\psi(q^a,q^p,q^{p-a})R^*(a,b,p;q)=P[q^A,q^B,q^{A+B}] 
\end{equation}
ii)
\begin{equation}
_2\phi_1[q^a,q^b;q^b;q^p,q^{(p-a-b)/2}]=\frac{\vartheta_4((a-b)i\log(q)/4,q^{p/2})}{\vartheta_4((a+b)i\log(q)/4,q^{p/2})}
\end{equation}
or
\begin{equation}
_2\phi_1\left[a,b;\sqrt{abc};c;\sqrt{\frac{c}{a b}}\right]=\frac{\vartheta_4((\log(a)-\log(b))i/4,c^{1/2})}{\vartheta_4((\log(a)+\log(b))i/4,c^{1/2})}
\end{equation}
\[
\]
\textbf{Proof.}\\
i)The proof of (61) follows easily from (58) and the q-binomial theorem \[
_2\phi_1[a,0,0,q,z]=\prod^{\infty}_{n=0}\frac{1-azq^n}{1-zq^n}
\]
ii)From the Gauss expansion
\begin{equation}
_2\phi_1[a,b;c;q,ab/c]=\frac{(c/a;q)_{\infty} (c/b;q)_{\infty}}{(c;q)_{\infty} (c/(a b);q)_{\infty}}
\end{equation}
we get the identity $_2\phi_1[q^{b-a},q^{a+b-p};q^b;q^p,q^{p-b}]=R^*(a,b,p;q)$. Hence
\begin{equation}
_2\phi_1[q^{b-a},q^{a+b-p};q^b;q^p,q^{p-b}]=\frac{\vartheta_4((p-2a)i\log(q)/4,q^{p/2})}{\vartheta_4((p-2b)i\log(q)/4,q^{p/2})}
\end{equation}
\[
\]
\textbf{Theorem 7.}\\
When $m_1$, $m_2$ are integres then 
\begin{equation}
R((2m_1+1)p/2,(2m_2+1)p/2,p;q)=1
\end{equation}
\begin{equation}
R(2m_1p,2m_2p,p;q)=1
\end{equation}
\begin{equation}
R(a,b,p;q)R(b,a,p;q)=1
\end{equation}
\[
\]
\textbf{Proof.}\\
It follows from Theorem 1 and the definitions of the Ramanujan`s Quantity.
\[
\] 
\textbf{Theorem 8.} 
\begin{equation}
R\left(-mp+\frac{i}{\sqrt{r}},p/2-mp+\frac{i}{\sqrt{r}},p;e^{-\pi \sqrt{r}}\right)=(-i)^m (k_{p^2r/4})^{1/2}
\end{equation}
where $r,p\in\bf R \rm$ and $m\in \bf N \rm$ 
\[
\]
\textbf{Proof.}\\
From Application 1 with transformation of Elliptic theta-4 function into Elliptic theta-3 and then using Theorem 1 and the definitions of the Ramanujan`s Quantity we get.
\begin{equation}
R(A,B,p,q^{2/p})=k_r^{1/2}(-i)^m
\end{equation}
where $A=-1/2(2m+c)p$, $B=-1/2(2m+c-1)p$, $c=i\pi/\log(q)$, $m\in N \bf$ and $p\in \bf R\rm$
the result follows from this equation.
\[
\]
\textbf{Examples.}\\
If $p$ real positive and $m$ integer, then
\begin{equation}
R\left(-p/2(2m+i),-p/2(2m+i-1),p;e^{-2\pi/p}\right)=(-i)^m 2^{-1/4}
\end{equation}
\begin{equation}
R\left(-(\sqrt{2}-4mi)pi/4,-(2-i\sqrt{2}-4m)p/4,p;e^{-\pi \sqrt{2}/p}\right)=(-i)^m \sqrt{\sqrt{2}-1} 
\end{equation}
\[
\]
\textbf{Corollary.}\\
If $\tau_0(a,q)$ is as in (41) and $m\in \bf N\rm$ then
\begin{equation}
\frac{\tau_0(m+1,q)}{\tau_0(m+1/2,q)}=(k_{r/4})^{1/2}
\end{equation}

\[
\]

\centerline{\bf References}\vskip .2in

[1]:I.J. Zucker, The summation of series of hyperbolic functions. SIAM J. Math. Ana.10.192(1979)

[2]:G.E.Andrews, Number Theory. Dover Publications, New York

[3]:T.Apostol, Introduction to Analytic Number Theory. Springer Verlang, New York

[4]:M.Abramowitz and I.A.Stegun, Handbook of Mathematical Functions. Dover Publications

[5]:B.C.Berndt, Ramanujan`s Notebooks Part I. Springer Verlang, New York (1985)

[6]:B.C.Berndt, Ramanujan`s Notebooks Part II. Springer Verlang, New York (1989)

[7]:B.C.Berndt, Ramanujan`s Notebooks Part III. Springer Verlang, New York (1991)

[8]:L.Lorentzen and H.Waadeland, Continued Fractions with Applications. Elsevier Science Publishers B.V., North Holland (1992) 

[9]:E.T.Whittaker and G.N.Watson, A course on Modern Analysis. Cambridge U.P. (1927)

\end{document}